\documentclass[12pt]{article}

\usepackage{graphicx} 
\usepackage{amsmath} 
\usepackage{amssymb,amsthm}
\usepackage{color}
\usepackage{tikz}
\usepackage{pgfplots}
\usepackage{enumitem}
\usepackage{mathtools}

\mathtoolsset{showonlyrefs=true}

\newcommand{\R}{{\mathbb R}}\newcommand{\N}{{\mathbb N}}

\newcommand{\Z}{{\mathbb Z}}\newcommand{\C}{{\mathbb C}}
\newcommand{\Q}{{\mathbb Q}}
\newcommand{\norm}[1]{\|#1\|}

\let\epsilon\varepsilon
\let\theta\vartheta
\let\hat\widehat

\newtheorem{theorem}{Theorem}[section]\newtheorem{lemma}[theorem]{Lemma}

\theoremstyle{definition}
\newtheorem{remark}[theorem]{Remark}

\title{About the occurrence of Akhmediev breathers in dispersive systems}
\author{Guido Schneider and Nils Thorin\footnote{Corresponding author: \texttt{nils.thorin@mathematik.uni-stuttgart.de}}\\
{\small
Institut f\"ur Analysis, Dynamik und Modellierung,} \\ {\small Universit\"at Stuttgart, Pfaffenwaldring 57, } \\ {\small 70569 Stuttgart, Germany}}

\begin{document}

\maketitle

\begin{abstract}
Akhmediev breathers are spatially periodic solutions of 
the NLS equation which are localized in time. 
We prove a new NLS approximation result for dispersive wave systems, showing that there are solutions to these systems where Akhmediev breathers in the lowest order approximately modulate the underlying wave packet.
Since the solutions of the original wave system are quasi-periodic in space, the main challenge is to find an appropriate function space to estimate the approximation error.
\end{abstract}

\section{Introduction}

The Nonlinear Schrödinger (NLS) equation 
\begin{equation}
    i \partial_T A + \nu_1 \partial_X^2 A + \nu_2 |A|^2 A = 0,
\end{equation}
with $ T, X \in \R $ and $ A(X,T) \in \C $ is a 
fundamental model in nonlinear physics. It describes
the evolution of complex wave envelopes 
in weakly nonlinear, dispersive media. It is a completely integrable
Hamiltonian system that can be solved explicitly using
an inverse scattering transform, cf. \cite{DJ89}. 
It possesses a number of special
solutions such as solitons, Peregrine solutions, Kuznetsov-Ma solitons,
or Akhmediev breathers. In previous publications it has been shown
that there are solutions of the original physical systems which behave 
as predicted by the NLS equation. Existing approximation theorems cover 
solitons, cf. \cite{Ka87,SU17book}, as well as recently, Peregrine solutions   and the 
Kuznetsov-Ma solitons \cite{ST25}.   However, they do not cover Akhmediev breathers. 
The goal of this paper is to close this gap.

\subsection{Akhmediev breathers}

\begin{figure}
  \centering
    \includegraphics[width=0.6\textwidth]{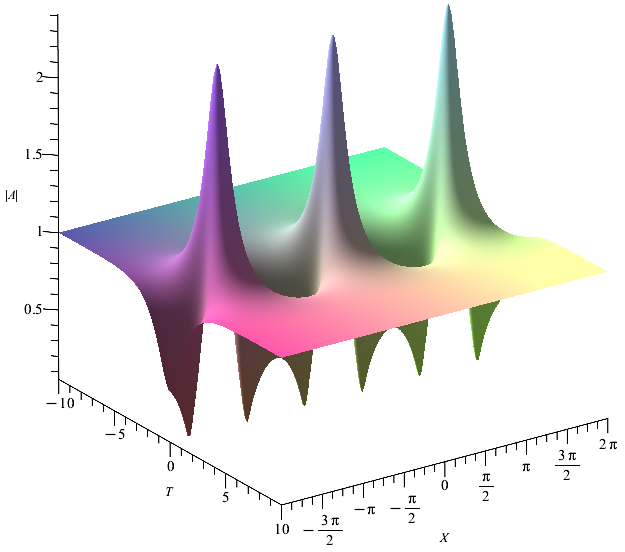}
    \caption{Absolute value of an Akhmediev breather  for $a=\frac{1}{4}$.\label{sec6fig1_a}}
 \end{figure}

For the 
focusing NLS equation in its normalized form 
$$
i \partial_\tau \psi+ \frac{1}{2}  \partial _\xi^2\psi+ |\psi|^2 \psi = 0
$$
there is a one-dimensional family of solutions
$$ 
\psi(\xi,\tau) = e^{i \tau}\left(1 + \frac{2(1-2a)\cosh(R \tau) + i R \sinh(R \tau)}{\sqrt{2a} \cos(L \xi)- \cosh(R \tau)} \right)
$$
where
$$ 
R = \sqrt{8a(1-2a)}, \qquad L = 2 \sqrt{1-2a}.
$$
For $ 0 < a < 1/2 $ the members of this family are called Akhmediev breathers.
They are spatially periodic and localized in time 
above a time-periodic, spatially constant background  state, cf. \cite{AAT09}.
See  Figure \ref{sec6fig1_a}.

\subsection{Quasiperiodicity}\label{sec:quasiperiodicity}
In this paper we will encounter so-called quasiperiodic (or almost-periodic) functions. We define a function $f:\R\to\R$ to be quasiperiodic if it can be written as (c.f. \cite{Sca08}) \begin{equation}
    f(x) = \sum_{m\in\Z^n}a_me^{i(\Lambda\cdot m)x},\quad \sum_{m\in\Z^n}|a_m|<\infty, \quad \Lambda\in\R^n.
\end{equation}
We especially note that if $\Lambda = (\lambda_1,\ldots,\lambda_n)$, then $f$ can be seen as the product of $n$ periodic functions, each with period $2\pi/\lambda_i$, $i=1\ldots n$. For our considerations, only $n=2$ will be used. 

Further, the only interesting case is when not all $\lambda_i$ are commensurable to each other, that is, it does not hold that $\lambda_i/\lambda_j\in\Q$  for all $i,j\in\{1,\ldots, n\}$. Otherwise, the quasiperiodic function degenerates into a regular periodic one.

\subsection{The NLS approximation}

The cubic Klein-Gordon equation
\begin{equation} \label{eq:3a1}
\partial_t^2 u = \partial_x^2 u - u + u^3,
\end{equation}
with $ x,t,u(x,t) \in \R $, is the simplest dispersive wave system where the NLS equation
can be derived.
The ansatz for the derivation is given by
\begin{equation} \label{ansatz1}
u(x,t) = \varepsilon \Psi_{NLS}(x,t) = \varepsilon A(\varepsilon(x-c_g t),\varepsilon^2 t) e^{i(k_0 x-\omega_0 t)} + c.c.,
\end{equation}
with spatial wave number $ k_0 > 0 $, temporal wave number $ \omega_0 \in \R $, group velocity $ c_g \in \R $, amplitude 
$ A(X,T) \in \C $, slow space and time variables $X=\varepsilon(x-c_g t)$ and $T=\varepsilon^2 t$, and 
small perturbation parameter $ 0 < \varepsilon \ll 1 $.
Inserting \eqref{ansatz1} into \eqref{eq:3a1} and equating the coefficients in front of $ \varepsilon^j e^{i(k_0 x-\omega_0 t)} $  to zero for $ j = 1,2,3 $ yields the linear dispersion relation $ \omega_0^2 = k_0^2 + 1 $, the linear group velocity $ c_g = k_0/\omega_0 $, and the NLS equation
\begin{equation} \label{nls5}
2 i \omega_0 \partial_{T} A = (c_g^2-1)\partial_{X}^2 A - 3 A |A|^2.
\end{equation}

\subsection{The approximation result}

The NLS equation was first derived by Zakharov in 1968 for the water wave problem \cite{Za68}.
Over the last few decades, NLS approximation results have been established for several dispersive wave systems.
In general, error estimates are nontrivial for quadratic nonlinearities, as solutions of order $\mathcal{O}(\varepsilon)$ must be controlled on a time scale of order $\mathcal{O}(1/\varepsilon^2)$.
The results are based on a near-identity change of variables that uses the oscillatory character of the  $\mathcal{O}(\varepsilon)$-terms in the error equations, cf. \cite{Ka87}.
For the above cubic Klein-Gordon equation 
an approximation result can be obtained by a simple application of Gronwall's inequality, cf. \cite{KSM92}.
However, for systems with quadratic nonlinearities resonant terms can lead to non-approximation results. There are rigorous counterexamples that the NLS approximation fails to make correct predictions, see \cite{SSZ15}.
All these approximation results hold for solutions which vanish for $ |x| \to \infty $. 
Recently, an approximation result \cite{ST25} was established for solutions of the form 
$ u = v + w $, where $v$ is spatially periodic and $w$ vanishes as $|x|\to\infty$. This result can be used to find  solutions to the cubic Klein-Gordon equation that are approximately described by the  Peregrine solutions and Kuznetsov-Ma solitons of the NLS equation. However, Akhmediev breathers are not covered by the existing approximation results.
Therefore, the goal of this paper is to prove the following theorem.
\begin{theorem} \label{thmain}
Let $A_{akh} $ be an Akhmediev breather 
of the NLS equation \eqref{nls5}.  Then for all $ T_0 > 0 $ there exist
$\varepsilon_0 > 0$ and $C>0$ such that for all $\varepsilon \in (0,\varepsilon_0)$
there are solutions $u$ of the original system \eqref{eq:3a1} which
can be approximated by $\varepsilon \Psi_{NLS} $ 
with
\[
\sup_{t\in [0,T_{0}/\varepsilon^2]} \sup_{x \in \R} |u(x,t)-\varepsilon \Psi_{NLS}(x,t)| < C\varepsilon^{2}.
\]
\end{theorem}
\begin{remark}{\rm 
For the NLS equation \eqref{nls5}  Akhmediev breathers $ A_{akh} $ are 
spatially periodic with a period $ 2\pi/L $, whilst the underlying carrier wave $ e^{i(k_0 x-\omega_0 t)} $ 
is spatially periodic with period $ 2 \pi/k_0 $. Hence, in general the NLS approximation $ \varepsilon \Psi_{NLS}(x,t) $
is quasi-periodic in space, and such solutions thus neither  vanish for $ |x| \to \infty $ nor  can be written as
$ u = v + w $ with $ v $ spatially periodic and 
$ w $ vanishing for $ |x| \to \infty $, and so the existing NLS approximation results do not apply.
Therefore,
the main challenge is to find an appropriate 
function space for estimating the error made by the approximation.
This will be done in Section \ref{sec2} followed by an adaption of the existing
NLS approximation proof for the cubic Klein-Gordon equation \eqref{eq:3a1}
to the chosen functional analytic set-up.
}
\end{remark}
\begin{remark}{\rm 
The presented approach can be combined with the existing NLS approximation theory which also includes systems with 
quadratic nonlinearities.   
Therefore, we expect similar approximation results  also for other nonlinear optics models and 
the water wave problem which will be the subject of future research.
In this paper, however, we restrict ourselves to the handling of the cubic Klein-Gordon equation since our focus is on the functional analytic set-up that can be used to prove an approximation result for Akhmediev breathers as solutions of the NLS equation.
}
\end{remark}
\begin{remark}{\rm 
Akhmediev breathers are essential for modeling modulation instabilities in optics and hydrodynamics. As a prototype for rogue waves, they explain the formation of extreme, localized energy concentrations in oceans and fiber optics \cite{CH24}.
}
\end{remark}
 \noindent
{\bf Notation.} Many, possibly different constants are denoted with the same symbol $ C $ 
if they can be chosen independently of the small perturbation parameter $ 0 < \varepsilon \ll 1 $.

 \noindent
{\bf Acknowledgement.}
The paper 
is partially supported by the Deutsche Forschungsgemeinschaft DFG 
through the SFB 1173 ``Wave phenomena" Project-ID 258734477.

\section{The function spaces}
\label{sec2}

We introduce the necessary function spaces for handling Akhmediev breathers of the NLS equation and their associated solutions of the cubic Klein-Gordon equation.

\subsection{The periodic function space}
As noted above, Akhmediev breathers are spatially periodic. Letting the period be $2\pi/L_0$, such functions may be written as a Fourier series of the form
$$ 
A(X,T) = \sum_{m \in \Z} \widehat{A}_m(T) e^{imL_0 X}.
$$ 
Here we collect some basic facts about these functions and introduce the relevant function space:
\begin{enumerate}[label=\textit{\roman*})]
    \item Multiplication in physical space corresponds to convolution in Fourier space. More precisely, given two periodic functions $A$ and $B$, we have \begin{equation}
        (\widehat{A B})_m  = (\widehat{A}*\widehat{B}  )_m = \sum_{n \in \Z} \widehat{A}_{m-n} \widehat{B}_n.
    \end{equation}
    
    \item A natural space to work in is the space of periodic function whose Fourier transform is in the space $\ell^1_s$, equipped with the norm \begin{equation}
     \| \widehat{A} \|_{\ell^1_s} = \sum_{m \in \Z} |\widehat{A}_m| (1+ |m|^s). 
    \end{equation}
    We further note that $\ell^1_s$ is closed under convolutions for $s\geq0$ due to the two inequalities \begin{align}
        \| \widehat{A}* \widehat{B}  \|_{\ell^1}  &\leq \| \widehat{A} \|_{\ell^1} \|  \widehat{B}  \|_{\ell^1}, \\
      (1+ |m+n|^s) &\leq C_s ((1+ |m|^s)+(1+ |n|^s)),
    \end{align}
    with the first one simply being  Young's inequality.
    
    \item The sup-norm of $A$ is naturally bounded by the $\ell^1_s$-norm of $\hat{A}$ for all $s\geq0$. Indeed: \begin{equation}\label{eq:sup_l1_ineq}
        \|A\|_{L^\infty} = 
    \sup_{X \in \R} |\sum_{m \in \Z} \widehat{A}_{m} e^{imL_0 X}| \leq \sum_{m \in \Z} |\widehat{A}_{m} | \leq \|\hat A\|_{\ell^1_s}.
    \end{equation} 
\end{enumerate}


\subsection{The quasi-periodic function space}

We recall that the NLS approximation $\varepsilon \Psi_{NLS}$ was written as \begin{equation}
    \varepsilon \Psi_{NLS}(x,t) = \varepsilon A(X,T) e^{i(k_0 x-\omega_0 t)} + c.c..
\end{equation}
Now, if the $A$ above is a spatially periodic Akhmediev breather with frequency $L_0$, then we see that our approximation is quasi-periodic with the two periods $2\pi/k_0$ and $2\pi/(\varepsilon L_0)$, where the extra $\varepsilon$ comes from the fact that $X = \varepsilon(x-ct)$.
Such functions can also be represented as a Fourier series 
$$ 
u(x,t) = \sum_{m,n \in \Z} \widetilde{u}_{mn}(t) e^{i(m \varepsilon L_0 +nk_0)x},
$$
and comparing to the definition given in Section \ref{sec:quasiperiodicity}, we see that $\Lambda = (\epsilon L_0, k_0)$.

To study functions of this form, we introduce the space 
$$
\ell^1_{p,q} = \{ \widetilde{u} :\mathbb{Z}^2 \to \C : \|  \widetilde{u} \|_{\ell^1_{p,q}} < \infty \} , 
$$
equipped with the norm 
$$ 
\|  \widetilde{u} \|_{\ell^1_{p,q}} = \sum_{m \in \Z} \sum_{n \in \Z}
| \widetilde{u}_{m,n} | (1+ |m|^p+|n|^q).
$$ 

In the following we prove that the NLS  approximation $ \varepsilon \Psi_{NLS} $
associated to 
Akhmediev breathers is $ \mathcal{O}(\varepsilon) $ in this space, cf. Lemma \ref{lemma1}, and
that this space is closed under convolutions, cf. Lemma \ref{lemma2}.
\begin{lemma} \label{lemma1}
Let $A_{akh} $ be an Akhmediev breather 
of the NLS equation \eqref{nls5}. 
Then for all $ T_0 > 0 $ and $ p,q \in \N_0 $ there exist
$\varepsilon_0 > 0$ and $C>0$ such that for all $\varepsilon \in (0,\varepsilon_0)$
we have for the associated NLS approximation that
$$
\sup_{t \in [0,T_0/\varepsilon^2]} \| \varepsilon  \widetilde{\Psi_{NLS}}(\cdot,t) \|_{\ell^1_{p,q}} \leq  C  \varepsilon.
$$
\end{lemma}
\noindent
{\bf Proof.} 
We represent the Akhmediev breather $ A_{akh}  $  of the NLS equation \eqref{nls5} with period $ 2\pi/L_0 $ as Fourier series
$$
A_{akh} (X,T) = \sum_{m \in \Z} \widehat{A}_m(T) e^{i m L_0X},
$$ 
with 
$ \sup_{T\in[0,T_0]}\sum_{m \in \Z} |\widehat{A}_m(T)| (1+ |m|^s) \leq C_s = \mathcal{O}(1) $
for all $ s \geq 0 $.
The associated NLS approximation  
$ \varepsilon \Psi_{NLS}(\cdot,t) $ is quasi-periodic  in space and can be represented as  
\begin{eqnarray*}
\varepsilon \Psi_{NLS}(x,t) & = & \sum_{m,n \in \Z} \varepsilon \widetilde{\Psi}_{m,n}(t) e^{imL_0 \epsilon x+ink_0x} \\
\\ & = & \varepsilon A(X,T) e^{ik_0 x + i \omega_0 t} + c.c. \\ & = & \sum_{m \in \Z} \varepsilon  \widehat{A}_m(T) e^{i m L_0 \varepsilon x-imL_0\varepsilon c_g t} e^{ik_0 x + i \omega_0 t} + c.c..
\end{eqnarray*}
Here we see after a comparison of coefficients, that \begin{align}
    \tilde\Psi_{m,n}(t) = \begin{cases}
        e^{i(\omega_0n-L_0c_g\varepsilon m)t} \hat A_m(\varepsilon^2t) ,\quad &\text{for } |n|=1,\\
        0, &\text{else}.
    \end{cases}
\end{align}
We can now estimate
\begin{eqnarray*}
\|  \varepsilon \widetilde{\Psi}_{NLS} \|_{\ell^1_{p,q}} & = &  \varepsilon \sum_{m,n \in \Z} |\widetilde{\Psi}_{m,n}(t) |  (1+ |m|^p+|n|^q) 
\\ & \leq & 2
 \varepsilon \sum_{m \in \Z} |\widehat{A}_m |  (1+ |m|^p+|1|^q)
 \leq 4 \varepsilon \| \widehat{A} \|_{\ell^1_p}.
\end{eqnarray*} 
\qed

Due to the cubic nonlinearity in the cubic Klein-Gordon equation \eqref{eq:3a1}, we will need to handle the multiplication of quasi-periodic functions. Taking two quasi-periodic functions
\begin{align}
u(x,t) &= \sum_{m,n \in \Z} \widetilde{u}_{mn}(t) e^{im \varepsilon L_0 x+ink_0x},\\ 
v(x,t) &= \sum_{m,n \in \Z} \widetilde{v}_{mn}(t) e^{im \varepsilon L_0 x+ink_0x},
\end{align}
and multiplying them gives 
\begin{eqnarray*}
u(x,t) v(x,t) & = & \sum_{m,n \in \Z}   \sum_{m',n' \in \Z} \widetilde{u}_{mn}(t) \widetilde{v}_{m'n'}(t)e^{i(m+m') \varepsilon L_0 x+i(n+ n')k_0x},
\end{eqnarray*} 
and so the multiplication of 
$ u $ and $v $ in physical space corresponds to the convolution
$$ 
(\widetilde{uv})_{m,n} =  (\widetilde{u} * \widetilde{v})_{m,n} = \sum_{m_1+m_2 = m}  \sum_{n_1+n_2 = n} \widetilde{u}_{m_1,n_1} \widetilde{v}_{m_2,n_2}
$$
in Fourier space. That this convolution is well defined in our chosen space is established in the following lemma.
\begin{lemma} \label{lemma2}
The space $ \ell^1_{p,q} $ is closed under convolution. In detail, 
for all $ p,q \in \N_0 $ there exists a $ C > 0 $ such that for all $ u,v \in  \ell^1_{p,q} $ we have 
$$ 
\| \widetilde{u} * \widetilde{v} \|_{\ell^1_{p,q}}  \leq C \| \widetilde{u} \|_{\ell^1_{p,q}} \|  \widetilde{v}  \|_{\ell^1_{p,q}} .
$$ 
\end{lemma}
\noindent
{\bf Proof.} 
The estimate follows from Young's inequality
$$ 
\| \widetilde{u} * \widetilde{v}  \|_{\ell^1}  \leq C \| \widetilde{u} \|_{\ell^1} \|  \widetilde{v}  \|_{\ell^1} 
$$
and the estimate 
$$ 
(1+ |m|^p+|n|^q) \leq C\left[(1+ |m-m'|^p+|n-n'|^q) + (1+ |m'|^p+|n'|^q)\right].
$$ 
\qed

\begin{remark}
{\rm For transferring the estimates from Fourier space to physical space we use, similar to \eqref{eq:sup_l1_ineq}, that for all $ p,q \in \N_0 $
\begin{equation} \label{fourierphysical}
\|u\|_{L^\infty} = \sup_{x \in \R} |\sum_{m,n \in \Z} \widetilde{u} _{m,n} e^{iml_0x+ink_0x}| \leq \sum_{m,n \in \Z} |\widetilde{u} _{m,n} | \leq \|\tilde u\|_{\ell^1_{p,q}}.
\end{equation}
}
\end{remark}

\begin{remark}
{\rm 
We also note that summability in Fourier space corresponds to regularity in physical space. In detail, from 
$$
    \partial_x^j u = \sum_{m,n\in\Z}i^j(m\epsilon L_0 + nk_0)^j\tilde{u}_{m,n}e^{i(m\epsilon L_0 + nk_0)x},
$$
we obtain 
$$
    \|\partial_x^j u\|_{L^\infty} \leq \sum_{m,n\in\Z} |\tilde u_{m,n}|(|m\epsilon L_0| + |nk_0|)^j \leq C_{j,L_0,k_0}\sum_{m,n\in\Z}|\tilde u_{m,n}|(1+|m|^j+|n|^j),
$$
where we used $ 0 < \epsilon \leq 1$ and $(|x|+|y|)^j \leq C_j(|x|^j+|y|^j)$. Therefore,  $\tilde u\in \ell^1_{p,q}$ implies that $u\in C^{\min(p,q)}$. 
}
\end{remark}
\begin{remark}{
The two parameters $p$ and $q$ allow us to control the decay on two scales, firstly on the $\mathcal{O}(1)$-scale via  $q$, and also on the $\mathcal{O}(\epsilon)$-scale via  $p$, see Figure \ref{fig:freq_distribution}. 
\begin{figure}
    \centering
    \includegraphics[width=0.7\linewidth]{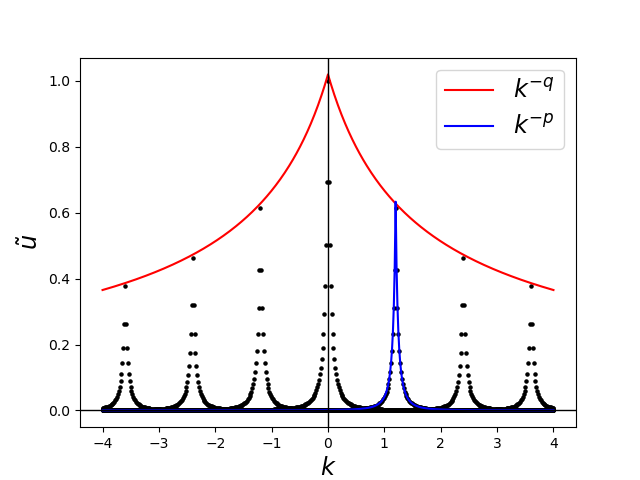}
    \caption{Frequency distribution of quasi-periodic functions in $\ell^1_{p,q}$}
    \label{fig:freq_distribution}
\end{figure}
}
\end{remark}

\section{The error estimates}

With the preparations from the previous section, we can now adapt the existing NLS approximation proof for the cubic Klein-Gordon equation 
as found in \cite[Section 11.2]{SU17book} for example. 

We introduce the error function $\varepsilon^{\beta} R$  with $\beta =  2$ by
$ u = \varepsilon \Psi + \varepsilon^{\beta} R $
where $ \varepsilon \Psi = \varepsilon \Psi_{NLS} + \mathcal{O}(\varepsilon^3) $
is an improved approximation subsequently defined in Lemma \ref{lemresi}. The choice of $\beta = 2$ will become apparent in the last part of this section. 
\begin{remark}\label{remresi}{\rm 
It is straightforward to see that Lemma \ref{lemma1} also holds for the improved 
approximation $ \varepsilon \Psi $. 
}
\end{remark}
In Fourier space the error function $ R $  satisfies
\begin{equation} \label{huerrv}
\partial_t^2 \widetilde{R} = -   \omega^2 \widetilde{R} +\varepsilon^2\widetilde{f} + \varepsilon^{- \beta} \widetilde{\textrm{Res}(\varepsilon \Psi)},
\end{equation}
where $ \omega =(\omega_{m,n})$, with $ \omega_{m,n} = \sqrt{(\varepsilon L_0m + k_0n)^2 + 1} $, 
\begin{eqnarray*}
\varepsilon^2 \widetilde{f}  & = &  3 \varepsilon^2 \widetilde{\Psi}^{*2}* \widetilde{R} + 3  \varepsilon^{1+\beta}
\widetilde{\Psi} *\widetilde{R}^{*2} + \varepsilon^{2 \beta} \widetilde{R}^{*3} ,
\end{eqnarray*}
and where the residual
$$ 
\textrm{Res}(u)  =  - \partial_t^2 u + \partial_x^2 u - u + u^3
$$ 
counts how much a function $ u $ fails to satisfy the cubic Klein-Gordon equation.
System \eqref{huerrv} is written as a first order system
\[
\partial_t \widetilde{\mathcal{R}}_{m,n}(t) = \widetilde{\Lambda}_{m,n} \widetilde{\mathcal{R}}_{m,n}(t) + \varepsilon^2 \widetilde{F}_{m,n}(t)+ \varepsilon^{- \beta} \widetilde{\textrm{RES}(\varepsilon \Psi)},
\]
with 
$$ 
\widetilde{\Lambda}_{m,n} = \left(\begin{array}{cc} 0 & i \omega_{m,n} \\  i \omega_{m,n} & 0 \end{array} \right), \qquad
 \widetilde{F}_{m,n}(t)= \left(\begin{array}{c} 0  \\  \frac{1}{i \omega_{m,n}} \widetilde{f}_{m,n}(t)\end{array} \right), $$ 
  and $$
\widetilde{\textrm{RES}(\varepsilon \Psi)} =  \left(\begin{array}{c} 0  \\  \frac{1}{i \omega_{m,n}} 
\widetilde{\textrm{Res}(\varepsilon \Psi)}\end{array} \right).
$$ 
We use the variation of constant formula
\begin{eqnarray*}
\widetilde{\mathcal{R}}_{m,n} (t) &=& e^{t \widetilde{\Lambda}_{m,n} } \widetilde{\mathcal{R}}_{m,n}(0) + \varepsilon^2 \int^t_0 
e^{(t-\tau)\widetilde{\Lambda}_{m,n}
  } \widetilde{F}_{m,n}(\tau) d \tau
 \\ &&  + \varepsilon^{- \beta} \int^t_0 
e^{(t-\tau)\widetilde{\Lambda}_{m,n}
  } \widetilde{\textrm{RES}(\varepsilon \Psi)}_{m,n}(\tau) d \tau.
  \end{eqnarray*}
to estimate the solutions of this system.
There are three ingredients to do so, namely estimates for the linear semigroup, cf. Lemma 
\ref{lemsemi}, estimates for the nonlinear terms, cf. Lemma 
\ref{lemnonl}, and estimates for the residual terms, cf. Lemma 
\ref{lemresi}.
%
%
%
%
%
%
%
%
\begin{lemma} \label{lemsemi}
The semigroup $(e^{t \widetilde{\Lambda}_{m,n}})_{t\geq 0}$ is uniformly bounded in every  $  \ell^1_{p,q} $, 
i.e., for every $ p,q \in \N_0  $
there exists a $ C > 0 $ such that for all $ \varepsilon \in (0,1) $ we have 
$$ \sup_{t\in \R} \|e^{t \widetilde{\Lambda}_{m,n}}\|_{\ell^1_{p,q} \to \ell^1_{p,q}} \leq C.$$
\end{lemma}
\noindent 
{\bf Proof.}
We have $\widetilde{\Lambda}_{m,n} = S \widetilde{D}_{m,n}S^{-1}$ and as a consequence
$
e^{t  \widetilde{\Lambda}_{m,n}} = S e^{t \widetilde{D}_{m,n}} S^{-1} $
where  
\[
S = \begin{pmatrix}1 & 1
    \\ 1 & -1 \end{pmatrix} \qquad {\rm  and} \qquad \widetilde{D}_{m,n} = \begin{pmatrix} i \omega_{m,n} & 0 \\ 0 & -i \omega_{m,n}\end{pmatrix}.
\]
The estimate then follows from
\[
\|e^{t \widetilde{\Lambda}_{m,n}} \widetilde{u}\|_{\ell^1_{p,q}} \leq \sup_{m,n\in\Z} \|e^{t \widetilde{\Lambda}_{m,n}
  }\|_{\C^2 \to \C^2} \|\widetilde{u}\|_{\ell^1_{p,q}}.
\]
\qed

\begin{lemma} \label{lemnonl}
For every $p,q \in \N_0 $ there is a $C > 0$ such that
for all $ \varepsilon \in (0,1] $ we have 
\begin{eqnarray*}
\|\widetilde{F}\|_{\ell^1_{p,q}} & \leq & C \left( \|\widetilde{\mathcal{R}}\|_{\ell^1_{p,q}} + \varepsilon^{\beta - 1}
  \|\widetilde{\mathcal{R}}\|^2_{\ell^1_{p,q}} + \varepsilon^{2\beta - 2} \|\widetilde{\mathcal{R}}\|^3_{\ell^1_{p,q}} \right).
\end{eqnarray*}
\end{lemma}
\begin{proof}
The estimate follows from 
$\|\frac{1}{\omega} \widetilde u\|_{\ell^1_{p,q}} 
\leq \|\widetilde u\|_{\ell^1_{p,q}} $, 
Lemma \ref{lemma1} together with Remark \ref{remresi}, and  Lemma \ref{lemma2}.
\end{proof}
\begin{lemma} \label{lemresi}
For the improved approximation 
\begin{eqnarray*}
u(x,t) = \varepsilon \Psi(x,t) & = & \varepsilon A(\varepsilon(x-c_g t),\varepsilon^2 t) e^{i(k_0 x-\omega_0 t)}  \\ && + 
 \varepsilon^3 A_3(\varepsilon(x-c_g t),\varepsilon^2 t) e^{3 i(k_0 x-\omega_0 t)}
+ c.c..
\end{eqnarray*}
with $ A $ satisfying the NLS equation \eqref{nls5} and $ A_3 $ determined by 
\begin{equation} \label{A3veq}
0 =   (9 k_0^2 - 9 \omega_0^2 + 1) A_{3} - A^3 
\end{equation} 
the following holds:
For  $ s_A -3 \geq p \geq 0 $, $ q \geq 0 $ and 
$ A  \in C([0,T_0],\mathcal{F}^{-1} \ell^1_{s_A}) $ 
there exist $ \varepsilon_0 > 0 $ and $ C > 0 $ such that 
for all $ \varepsilon \in (0,\varepsilon_0) $ we have 
$$ 
\sup_{t \in [0,T_0/\varepsilon^2]} \|\widetilde{\textrm{\rm Res}( \varepsilon \Psi) }\|_{\ell^1_{p,q}}
\leq C \varepsilon^4 .
$$ 
\end{lemma}
\begin{proof}
Inserting the improved approximation into \eqref{eq:3a1}
and equating the coefficient in front of 
$  \varepsilon^j e^{i(k_0 x-\omega_0 t)} $ for $ j =1,2,3 $
and
$  \varepsilon^3 e^{3 i(k_0 x-\omega_0 t)} $ to zero gives the formulas 
from the introduction, namely the  
dispersion relation, a formula for the group velocity,  
NLS 
equation \eqref{nls5} and finally the relation \eqref{A3veq}.
By this choice all terms of order $ \mathcal{O}(\varepsilon^3) $ in the residual cancel such that 
it is formally of order  $ \mathcal{O}(\varepsilon^4) $.
For making these formal estimates rigorous we adapt \cite[Lemma 11.5.1]{SU17book}
to $ \ell^1_{p,q} $-spaces.

We also take inspiration from \cite{LS26}, by noting that functions in $\ell^1_{p,q}$ can be expressed as  functions on $\R$ in the sense that \begin{equation}
    \norm{\tilde u}_{\ell^1_{p,q}} = \sum_{m,n\in\Z}|\tilde u_{m,n}|(1+|m|^p+|n|^q) = \int_{\R} |\tilde u(k)|d\mu_{p,q}(k),
\end{equation}
with $d\mu_{p,q}(k) = \sum_{k'\in\Lambda}(1 + |(L_0\epsilon)^{-1}(k'\cdot e_1)|^p + |k_0^{-1}(k'\cdot e_2)|^q)\delta(k-(k'\cdot\mathbf{1}))$ being a discrete measure with weighted point masses $ 1 $ on the points $\{k_0n+\epsilon L_0m;\, m,n\in\Z\}$, with $\Lambda$ being  the lattice $\Lambda = \{(k_0n,\varepsilon L_0m);\;(m,n)\in\Z^2\}$.

This means that we may evaluate functions defined on $\R$ in $\ell^1_{p,q}$ by simply looking at the values of the function at the points $k_0n+\epsilon L_0 m$. This allows us to transfer the known results fairly painlessly.

\begin{lemma} \label{estilemdd} 
Let $ p,q, {\theta}_0\geq 0 $, $ {\theta}_{\infty} \in \R $, and let 
$ g(k) $ satisfy 
$$ | g(k) | \leq C \min (| k -
k_0 |^{{\theta}_0}, (1+|k|)^{{\theta}_{\infty}}).$$
Then 
$$ s = \| g(\cdot)\epsilon^{-1}\widehat{A}(\epsilon^{-1} (\cdot - k_0)) \|_{\ell^1_{p,q}}
\leq C \varepsilon^{{\theta}_0} \| \widehat{A} \|_{\ell^1_{p+\max({\theta}_{0},{\theta}_{\infty})}} .$$  
\end{lemma}
\begin{proof}  We have 
\begin{eqnarray*}
s & = &\sum_{m \in \Z} 
| g(k_0 + \varepsilon L_0 m) \widehat{A}_m | (1+ |m|^p+|1|^q)
\\ 
& \leq & C \sum_{m \in \Z} 
 \min (| \varepsilon L_0m |^{{\theta}_0}, (1+|k_0 + \varepsilon L_0m|)^{{\theta}_{\infty}}) |\widehat{A}_m | (1+ |m|^p+|1|^q)
\\
& \leq & C \varepsilon^{{\theta}_0} \|  \widehat{A} \|_{\ell^1_{p +\max({\theta}_{\infty},{\theta}_0)}}
\end{eqnarray*}
\end{proof}

With this lemma the formal expansion of the curve of eigenvalues and of the kernels 
in the multilinear maps can be estimated, where given that $\omega(k) = \sqrt{1+k^2}$ we let \begin{equation}
    g(k) = \omega(k) - \omega_0 - c_g(k-k_0) - \frac{\omega''(k_0)}{2}(k-k_0)^2 = \mathcal{O}(|k-k_0|^3).
\end{equation} The rest of the proof of Lemma \ref{lemresi} is then 
 completely straightforward, cf. \cite[Section 11.2]{SU17book}. 
\end{proof}

Now, using Lemma \ref{lemsemi}, Lemma \ref{lemnonl},  and Lemma 
\ref{lemresi} shows that
the quantity 
$ 
Z(t) = \|\widetilde{\mathcal{R}}(t)\|_{\ell^1_{p,q}} 
$
satisfies the inequality
\begin{align*}
Z(t)
&\leq CZ(0) + C\varepsilon^2 \int^t_0 \left( Z(\tau)+
  \varepsilon^{\beta- 1} Z(\tau)^2 + \varepsilon^{2\beta-2}
 Z(\tau)^3 + 1\right) d \tau\\
&\leq  CZ(0)+ C\varepsilon^2 \int^t_0 \left( Z(\tau) +
  2 \right) d \tau
\leq C(\tilde{C}+2T_0) + C \varepsilon^2 \int^t_0 Z(\tau) d
\tau  
\end{align*}
which holds as long as 
\begin{equation} \label{witt}
\varepsilon^{\beta -1} Z(\tau)^2+
\varepsilon^{2\beta- 2}Z(\tau)^3 \leq 1,
\end{equation}
 where
$\tilde{C} = Z(0) $.
Applying Gronwall's inequality yields
\[
Z(t) \leq C(\tilde{C}+2T_0) e^{C\varepsilon^2t} \leq C(\tilde{C}+2T_0)
 e^{CT_0} = M
\]
for all $t\in [0,T_{0}/\varepsilon^2]$. Choosing $\varepsilon_0 > 0$ such that
$\varepsilon_0^{\beta - 1} M^2 + \varepsilon_0^{2\beta - 2} M^3 \leq 1$ 
ensures that condition \eqref{witt} is satisfied. 
An application of \eqref{fourierphysical}
completes the proof of our approximation result.

\begin{remark}{\rm 
The choice of $\beta = 2$ can now be motivated in hindsight, 
since the entire integral above must be of the order of $\mathcal{O}(\epsilon^2)$ in order to achieve the estimates on the long time scale $\mathcal{O}(\epsilon^{-2})$.
}
\end{remark}

\section{Discussion}

With some obvious modifications to the proof in \cite{ST25}, we can consider solutions of the form $ u = v+ w $,
where $ v $ is quasi-periodic in space and $ w $ vanishes as $ |x| \to \infty $. 
We write the cubic Klein-Gordon equation for $v$ and $w$ as  
\begin{eqnarray} \label{veq}
\partial_t^2 v & = & \partial_x^2 v - v + v^3 ,
\\
\partial_t^2 w & = & \partial_x^2 w - w + 3 v^2 w + 3 v w^2 + w^3 .\label{weq}
\end{eqnarray}
The first equation has been discussed in previous sections. Since the product of a quasi-periodic function $v$ and a localised function $ w $ is still a localised function, as in  \cite{ST25}, we can conclude
\begin{theorem} \label{thmain2}
Let $ p,q \geq 0 $, $s \geq \max(p+ 3,q) $ and
let $ A = A_v + A_w$, with 
$$(A_v, A_w) \in C([0,T_0],\mathcal{F}^{-1} \ell^1_{s+3} \times H^{s+3}(\R)) ,$$ 
be a solution of the NLS equation \eqref{nls5}.
Then there exist $ \varepsilon_0 > 0 $ and $ C > 0 $ such that for all
$ \varepsilon \in (0, \varepsilon_0) $ there are solutions $ u = v+ w $ with 
$$(v,w) \in C([0,T_0/ \varepsilon^2],\mathcal{F}^{-1} \ell^1_{p,q}\times H^s(\R))$$ 
of the cubic Klein-Gordon  equation \eqref{eq:3a1} satisfying
$$
\sup_{t \in [0,T_0/ \varepsilon^2]} \| (v,w)(\cdot,t) - (\varepsilon\Psi_v,\varepsilon\Psi_w)(\cdot,t) \|_{\mathcal{F}^{-1} \ell^1_{p,q} \times H^s} \leq C \varepsilon^{3/2},
$$
where $ \varepsilon\Psi_v $ and $\varepsilon\Psi_w $ are the associated NLS approximations of 
$ v $ and $ w $, respectively.
\end{theorem}
We expect the existing theory on NLS approximations, such as the handling of quadratic nonlinearities via normal form transformations (cf. \cite{Ka87}), or of resonances (cf. \cite{Schn05}), or of quasilinear systems (cf. \cite{Du17,WC17,DH18}), to be transferable.
For example, we expect the approximation result for the water wave problem to be similar to that in Theorem \ref{thmain}. 
For the handling resonances or quasilinear systems, the $\ell^1_{p,q}$-spaces 
still have to be used for the approximation, but for the error, we have to use 
$\ell^2_{p,q}$-spaces. This will be subject of future research.

\bibliographystyle{alpha}
\bibliography{perebib}

\end{document}